\documentclass{article}
\usepackage[utf8]{inputenc}

\title{Invariance under permutations as a semantic motivation for Stratification}
\author{Zuhair Al-Johar}
\date{24 October 2020}

\begin{document}

\maketitle

\section{Abstract}
This article examines the notion of invariance under different kinds of permutations in a milieu of a theory of classes and sets, as a semantic motivation for Quine's new foundations ``NF". The approach largely depends on interpreting a finite axiomatization of NF beginning from the least restrictions on permutations and then gradually upgrading those restrictions as to enable interpreting NF. Comparisons with earlier works are drawn.

\tableofcontents

\section{Introduction}

That stratified formulas are those found invariant under setlike permutations of the universe is a topic that had been extensively studied by Thomas Forster[6]. Randall Holmes building on Forster's work began the project of interpreting stratified comprehension in a theory about symmetric sets [7] [8]. In the last of his articles he exactly phrased that in an NBG style class\textbackslash set theory. This article copies the general approach of Holmes's last article, but it defines matters in a different way. The proof of interpreting Stratified Comprehension given here is through proving a finite axiomatization of stratified comprehension that I've previously presented [1]. The notion of invariance given here might provide a semantic insight into stratified theories, yet much of the development here is still technical.

 \section{Main theory}
 
 This theory is written in mono-sorted first order logic with equality $=$ and class membership $\in$. \bigskip

\emph{Define:} $set (x) \iff \exists k (x \in k)$ \bigskip

\noindent
1. \textbf{Extensionality:} $\forall z (z \in X \leftrightarrow z \in Y) \rightarrow X=Y$\smallskip

\noindent
2. \textbf{Class Comprehension:} if $\phi$ is a formula in which $x$ is not free, then $(\exists x \forall y (y \in x \leftrightarrow set(y) \land \phi))$ is an axiom.\smallskip

\noindent
\emph{Define} $V$ as the class of all sets.\smallskip

\noindent
Define: $x = \{y \in V: \phi\} \iff \forall y (y \in x \leftrightarrow set(y) \land \phi)$ \smallskip

\noindent
3. \textbf{Pairing:} $\forall a \forall b \forall x: \forall y \in x (y=a \lor y=b) \to set(x)$\bigskip

\noindent
Define $functions$ in terms of classes in the usual manner as classes of Kuratowski ordered pairs (denoted as $\langle , \rangle$), wherein every two pairs having a common first projection must have the same second projection. \bigskip

\noindent
\emph{Define:} $f \ permute \ X \iff f: X \leftrightarrow X $

Where $``\leftrightarrow"$, denotes ``is a bijective function from to" \smallskip

\bigskip

\noindent
\emph{\textbf{Define($f$-membership):}} $$y \in^f  x \iff  f(y) \in x$$, where $f$ is a function. \bigskip

A formula $\phi^f$ is meant to be the formula obtained from the set theoretic formula $\phi$ by merely replacing each occurrence of the membership symbol $\in$ by the $f$-membership symbol $\in^f$.  \bigskip

\noindent
We spell out the invariance notion, this method is all about: \smallskip 

A formula $\phi$ having all of its free variables among symbols $``y,x_1,..,x_n"$, is said to be $invariant$ under $f$, if and only if: $$\forall x_1,..,x_n \forall y (\phi \leftrightarrow \phi^f)$$
We first examine the most naive form of invariance notion, which I'll spell in terms of the following schema. The idea is that if a class $x$ is defined by an instance of naive comprehension, that is as: $\forall y (y \in x \iff \varphi)$ with parameters from $\vec{p}$, take this formula to be $\phi$, and the definition was invariant under all permutations $f$ of $V$, then $x$ is a set (i.e. $x \in V$). \bigskip

\noindent
\textbf{\emph{Naive invariance schema}}: For the above qualifications: $$\forall \vec{p} \in V \forall x (\phi \land \forall f ([f: V \leftrightarrow V]  \rightarrow [\phi \leftrightarrow \phi^f]) \Rightarrow x \in V)$$ Now of a finite axiomatization of SF, this schema would prove the following sentences:

$V \in V$

$\forall a \in V :a^{\mathsf c} \in V$

$\forall a,b \in V: a \cup b \in V $ \bigskip

\noindent

Where $a^{\mathsf c}$ is the absolute complementary set of $a$. This schema is weak for the purpose of interpreting $SF$. So we need to employ some restriction on the kind of permutations the invariance under which results in sets. This would be the intention behind the following development! The general idea is that we'll require the permutations to fix several stages of the set structure, and examine invariance under those conditions. \bigskip

Now when $f$ permutes a class $x$, then it would fix its $f$-membership (i.e.; its fixes $x$ class-wisely), that is its $f$-members would be exactly its $members$! $$ Define: f``x=\{f(y): y \in x\} $$, so when $f$  permute $x$, then  we'll have $f``x \subset x$ and $ x \subset f``x$, so $f``x=x$, and so $x=\{y: y \in^f x\}$, also $x=\{f(y): y \in^f x\}$. accordingly: $\forall y (y \in x \iff y \in^f x)$\bigskip

Define: $\{x:\phi\}^f \iff \{f(x):\phi\}$ \bigskip

so if $f \ permute \ x$, then we'll have: $x=\{y: y \in^f x\}^f$ \bigskip

For any permutation $f$ of $V$ we define: \bigskip

$j`f = \{\langle x, f``x \rangle: x \in V\}$ \smallskip

$j^0`f=f; j^{n+1}`f = j`j^n`f $ \bigskip

\newpage

\emph{Define:}  

$f \ 0$-$permute \ x \iff f: V \leftrightarrow V$

$f \ 1$-$permute \ x \iff f \ 0$-$permute \ x \land  f \ permute \ x $

$f\ 2$-$permute \ x \iff f \ 1$-$ permute \ x \land j^1`f  \ permute \ x$ \smallskip

$2$-permutation shall be denoted as $``ultrapermutation"$. \bigskip

\noindent
Our intention is to employ invariance under all $f$ permutations of the kinds defined above, so we'll employ those permutations over the set naively defined and over all parameters used in that definition as well. This would enable us to prove more theorems, and it'll be seen to be strong enough to interpret all axioms of NF.\bigskip

\noindent
4. \textbf{Axiom schema of invariant comprehension:} If $\phi$ is the formula $``\forall y (y \in x \leftrightarrow y \in V \land \varphi)"$ where $\varphi$ is a formula in the language of set theory with all quantifiers bounded $\in V$, and in which only variable symbols $y,w,v$ occur free,  let $\phi^f$ be the formula obtained from $\phi$ by merely replacing each symbol $\in$ in it by $``\in^f"$; then the following is an axiom: \bigskip

\noindent
$\ i=0,1; j=0,2 ; \forall w, u \in V; \forall x \\\big{[} \phi \land \forall f ( f \ i$-$ permute \ w,v \land  f \  j$-$permute \ x \rightarrow [\phi \leftrightarrow  \phi^f]) \Rightarrow   x \in V  \big{]}$ \bigskip

/ Theory definition finished.

\section {Interpreting NF}

\noindent
The method is to prove the five axioms in my latest finite axiomatization of $SF$ [1].\smallskip

I'd first spell out the 5 axioms of Stratified Comprehension: \smallskip

1. Sheffer stroke: $\forall a,b \in V \exists c \in V: c= \ a \uparrow b $

2. Singletons: $\forall a \in V \exists b \in V:  b=\{a\} $

3. Set Union: $\forall a \in V \exists b \in V: b= \bigcup(a) $

4. Unordered relative products: 

$\forall r,s \in V \exists x \in V: x= \{\{a,c\}: \{a,b\} \in r \land \{b,c\} \in s \}$

5. Unordered intersection relation set: 

$\exists k \in V: k= \{\{x,y\} :\exists z (z \in x \land z \in y)\}$\bigskip

\noindent
Singletons is already an axiom of our system. The Sheffer stroke set can be proven from Boolean union and Complements already theorems of naive invariance (which results from fixing $i,j$ to 0 in the above scheme). \smallskip

\noindent
Set unions begs $1$-permutation of a parameter by the $f$s' and that's enough to ensure invariance under all such permutations, because by 1-permutation of $A$ we have: \bigskip

\noindent
$\bigcup(A) = \{x: \exists y (y \in A \land x \in y)\} = \{x: \exists y (y \in^f A \land x \in y)\}$. \bigskip

\noindent
Now since $f \ permute \ V$, then we'll have: \bigskip

\noindent
$\{f(z): \exists y (y \in^f A \land f(z) \in y)\}= \{z: \exists y (y \in^f A \land z \in^f y)\}^f $\bigskip

\noindent
To get the unordered relative products we need $1$-permutation level over parameters and the defined set as well. The reason is because we can arrange for $f$ to permute $V$ and at the same time permute any sets $R,S,K$. Now if $K$ is the unordered relative product of $R,S$, i.e., $K = R:S$; then we have:$$\forall y (y \in K \iff \exists a,b,c (\{a,b\} \in R, \{b,c\} \in S, y= \{a,c\}))$$;  since $f$ is permuting $R,S, K$; then:$$\forall y (y \in^f K \iff \exists a,b,c (\{a,b\} \in^f R, \{b,c\} \in^f S, y= \{a,c\}))$$; since $f$ also permutes $V$, then we have: \bigskip

\noindent
$\exists a,b,c (\{a,b\} \in R, \{b,c\} \in S,  y= \{a,c\}) \iff \\\exists a,b,c (\{f(a),f(b)\} \in R, \{f(b),f(c)\} \in S,  y= \{f(a),f(c)\}) $. \smallskip

Then we have:  $$\forall y (y \in^f K \iff \exists a,b,c (\{a,b\}^f \in^f R, \{b,c\}^f \in^f S, y= \{a,c\}^f))$$ 
Ultrapermutation is only needed to prove the unordered intersection relation set, so here there would be no parameters. Define the unordered-intersection relation class $\{\cap\}$ as the class of all unordered pairs of intersecting elements, that is: $\{\cap\}=\{\{ x,y \}: \exists c (c \in x \land c \in y)\}$. Now if $f$ ultrapermute $\{\cap\}$ then all elements of $\{\cap\}$  would be $f$-elements of $\{\cap\}$   [by 1-permutability of $f$], now take any $\{A,B\} \in \{\cap\}$, we'll have $\{f^{-1}(A),f^{-1}(B)\} \in \{\cap\}$ [by ultrapermutability of $f$ over $\{\cap\}$], thus we'll have $(\exists x: x \in f^{-1}(A) \land x \in f^{-1}(B))$ [definition of $\{\cap\}$], Since $f$ permute $V$, then for any shared $\in$-element $x$ of $f^{-1}(A)$ and $f^{-1}(B)$, there will exist a set $y$ such that $x=f(y)$, by then we'll have $y \in^f f^{-1}(A) \land y \in^f f^{-1}(B)$ [Definition of $\in^f$], so any element of $\{\cap\}$ is an unordered pair whose $f$-elements are $f$-intersecting, establishing the invariance over $\{\cap\}$, thus $\{\cap\} \in V$. An example of a non-trivial ultrapermutation over $\{\cap\}$ is any automorphism over $V$, and $\{\cap\}$ would be preserved (invariant) under all of those. \bigskip

So all axioms of SF are provable in this theory. We already have Extensionality, then we have NF. \bigskip

\noindent
I don't know what's the consistency status of this theory relative to NF or to Randall Holmes's last theory on symmetric comprehension[7].

\section{Avoiding Known paradoxes}

1. Russell's paradox: The class $R$ of all sets that are not members of themselves, is not a set. Formally: $$\{x \in V: x \not \in x\} \not \in V$$
From Class Comprehension we do have the class $R$ of all sets that are not members of themselves (i.e. the Russell class) existing! We'll show that this is NOT invariant under all 0,1-permutations of it, and under all ultra-permutations of it.

Let $f$ send $\emptyset$ to $\{\emptyset\}$ and $\{\emptyset\}$ to $\emptyset$, and fix all sets otherwise. Then $f$ fixes $R$ class-wisely, i.e. $f``R=R$, yet $R$ doesn't satisfy $\forall y (y \in^f R \Leftrightarrow y \not \in^f y)$ since $\{\emptyset\} \in^f R$ and yet $\{\emptyset\} \in^f \{\emptyset\}$! Now $f$ is both 0-permutation and 1-permutation over $R$. To check ultrapermutations, take an automorphism $\pi$ over $V$, this will preserve the membership structure of $R$, so it'll be an ultrapermutation over $R$, now let $\pi$ send an ordinal to a preceding one (which exists in some model of this theory [5] ), then this ordinal would be  $\pi$-element of itself! \bigskip

2. Cantor's paradox: $\forall x \in V (|P(x)| > |x| \Rightarrow V \not \in V)$

Proof: The premise of the above statement is not provable in this theory because it depends on $Separation \ for \ Sets$ (i.e.; any definable subclass of a set, is a set), which is not true in this theory! So we can have $V \in V$, and of course $P(V)=V$, so the identity function I on $V$ is a bijection between $V$ and $P(V)$, and the diagonal class $\{x \in V : x \not \in I(x)\}$ is exactly the Russell class, but this is NOT a set! \bigskip

3. Leśniewski's paradox: The class of all singletons that are not elements of their sole members cannot be a set. Formally: $$S=\{\{x\} \in V:  \{x\} \not \in x\} \not \in V$$ 

Let $f$ send $\{\emptyset\}$ to $\{\{\{\emptyset\}\}\}$ and vice verse, and fix all other sets. So $f$ is a 0,1,ultra-permutation over $S$. Now clearly we'll be having $\{\{\{\emptyset\}\}\} \in^f S \land \{\{\emptyset\}\} \in^f \{\{\{\emptyset\}\}\} \land \{\{\{\emptyset\}\}\} \in^f \{\{\emptyset\}\}$. !\bigskip

4. Mirimanoff paradox: The class of all well founded sets cannot be a set. Proof: the first function depicted above for Russell's paradox will be a 0,1,ultra-permutation, and the same argument runs as for Russell's paradox. \bigskip

5. Burali Forti paradox: Two versions, the Von Neumann ordinals version, and the Frege ordinals version.

5.a: The von Neumann ordinals version: The class of all von Neumann ordinals that are sets, is not a set. Formally: $$ON=\{\alpha \in V: von \ Neumann (\alpha)\} \not \in V$$ Proof: take an automorphism $f$ over the universe of this theory that moves some ordinals inwards, i.e.; for some $\alpha: f(\alpha) < \alpha$, it's known that such automorphism must exist for some model of the theory. Now $f$ is a (0,1,ultra)-permutation over $ON$, and clearly all moved elements of $ON$ (non-standards) would be $f$-members of themselves, and all of those are $f$ members of $ON$, therefore $ON$ is no longer a class of $f$-ordinals. Now this happened because $ON$ contains non-standard ordinals. However, if we define $ON$ as the class of all TRUE ordinals, i.e. those that are truly well founded by $\in$, then no non-trivial permutation can exist over it, but this class is not eligible for the invariance scheme since its definition encounters quantification over non-sets. A proof of the same result without resorting to automorphisms is also easy: take the function presented in the proof of Russell's paradox, this is $0,1$-permutation over $ON$, that is not invariant over $ON$ since we have $\{\emptyset\} \in^f ON$ and also $\{\emptyset\} \in^f \{\emptyset\}$. Now for any non-trivial ultra-permutation over $ON$ which can only happen because there are non-standards in $ON$, the set of all of which doesn't have a least element with respect to membership $\in$, now if $f$ moves one ordinal inwardly then the moved ordinal would be an $f$-element of itself; now if $f$ only moves ordinals upwardly, then its converse would also be an ultra-permutation over $ON$ that moves some ordinals inwardly.

 5.b: The Frege ordinals version: the order type of the set of all order types that are strictly smaller than the order type $\Omega$ of the class $ON$ of all order types that are sets (equivalence sets of set-well orderings under order-isomorphism), is strictly smaller than $\Omega$. 

This is provable in NF, which is interpretable in this theory. Burali-Forti is prevented because the class $ON$ is not truly well ordered! That is there is no true well ordering over it. The definition of a true well ordering is prevented in this theory through restricting quantifiers and parameters in the definition of sets to sets only, thereby precluding Rosser's argument for Burali-Forti in the context of ML [11], in the same manner as how the corrected version of Quine's ML [10] prevents it as shown by Hao Wang [12]. 

\section{Avoiding challenges to symmetry}

Randall Holmes had presented two main challenges to the notion of symmetry, and here this theory shares the same challenges:

1. Paradox of the class of all ordinals of true well orderings: The class of all Russell-Whitehead ordinals of true well orderings (well orderings having every subclass (not just every subset) of their range possessing a minimal element [13]) is indeed invariant under all (0,1,2)-permutations. But its not eligible to be a set in this theory through the invariance axioms because its definition requires quantification over classes.

2. Paradox of strongly Cantorian sets.
Holmes had posed the question about the class STC of strongly Cantorian sets [in relation to SST], pointing the problem that STC is invariant under $set$ permutation, so this might constitute a challenge to this method as well. However, this is resolved here: \smallskip

Proof: first we define what does it mean for a class $x$ to be strongly cantorian set, so we define: $$x \ is \ a \ strongly \ cantorian \ set \iff x \in V \land  \exists f \in V f=\{\langle y, \{y\} \rangle: y \in x\})$$, the idea is to use class automorphisms, so we take a model $M$ of our theory, take an automorphism $f$ on it, also its known to preserve the strongly cantorian property, so it'll fix the class STC of all strongly Cantorian sets class-wisely (i.e. it'll permute STC), also by automorphism it'll be an ultra-permutation over STC. Now take $f$ to be a rank shifting automorphism, i.e. $f$ moves all the ranked sets in $M$ inwardly. Let $V_\alpha$ be some limit rank in $M$ that is in STC, so it'll be moved by $f$ inwardly to a lower limit rank $V_{f(\alpha)}$ which $M$ itself would see as *strictly smaller* than $V_\alpha$. Accordingly from the viewpoint of $M$ we'll have: $$f(\iota`` V_{f(\alpha)}) < (\iota`` V_{f(\alpha)}) < \iota`` V_\alpha$$, Where: $\iota``x = \{\{y\}: y \in x\}$ \bigskip

So M would see:$$f(\iota`` V_{f(\alpha)}) <  V_{f(\alpha)}$$ So we have: $$\{ \{ f(x) , f(\{f(x)\}) \} : x \in V_\alpha \} \not \in M $$ So by automorphism, we have:

$$\{ f( \{ f(x) , f(\{ f(x) \}) \} ) : x \in V_\alpha \} \not \in M$$ And so $V_\alpha$ is NOT an  $f$-strongly Cantorian $set$! \bigskip

\noindent
So STC is not invariant under ultra-permutation! QED \bigskip

if we remove the restriction on $f$ to be a set in the definition of strongly cantorian sets, then this proof won't work, but by then the defined class would not be eligible for being used in the invariance scheme.

\section {Comparison with earlier works}

The idea of equating stratification with invariance under certain kind of permutations thereby giving a semantic equivalent to the usual syntactical approach of stratification is first given by Forster 1995. 
That work was then followed with works on symmetric set comprehension by Randall Holmes. However, the invariance notion presented here is different from the one presented there. There $\phi^f$ would be defined after pre-fixing variables in $\phi$ after the $j^n`f$ functions, in such a manner that for every formula $`x \in y$' in $\phi$ if $`x$' is prefixed as $j^n`f(x)$, then $`y$' would be prefixed as $j^{n+1}`f(y)$, while for a formula $`x = y$' both occurrences of $`x$';$`y$' in it receive the same $j^n`f$ pre-fixing. The value of $n$ in the prefix is determined by the stratification given to the variable in a formula (or more generally by a kind of indexing on occurrences of variables in any formula that would yield a stratification should the formula be stratified) So stratified formulas could be identified with those that are invariant (i.e. $\phi \longleftrightarrow \phi^f$) under all $setlike$ permutation of the universe $V$. For that principle to work as intended, its important that permutations must be restricted to $setlike$ ones, where a permutation $\pi$ of $V$ is $setlike$ if and only if $j^n`\pi \ permute \ V$ for every meta-theoretic natural $n$.\smallskip

Invariance in the current method is not restricted to $setlike$ permutations. It requires certain class fixings that permutations must exert on the parameters and the defined classes for the invariance notion to work under, which are generally different from the earlier approaches. So the methods are not the same. A version of this method that would share some features of their method is to state invariance (in our style) under all automorphisms that fix parameters. Formally this is: \bigskip

\noindent
$\forall a,b \in V, \forall x: \\\phi(x) \land \forall f (auto(f) \land f(a)=a \land f(b)=b \rightarrow[\phi(x) \iff \phi^f(x)] ) \Rightarrow x \in V$ \bigskip

where $\phi(x)$ is the formula $``\forall y (y \in x \leftrightarrow y \in V \land \varphi(y,a,b))"$, with all quantifiers in $\varphi$ bounded $\in V$; and $$ auto(f) \iff f: V \longleftrightarrow V \land \forall x,y \in V (x \in y \iff f(x) \in f(y))$$ This would avoid all the known paradoxes listed above, the paradox of singletons would be avoided by using a rank shifting automorphism $f$ that moves ranks inwardly, take any stage $V_\alpha$ where $\alpha$ is some nonstandard limit ordinal, now $\{ V_\alpha\}$ would be sent by $f$ to $\{ V_{f(\alpha)}\}$ which is an element of $V_\alpha$, so $\{ V_\alpha\} \in^f V_\alpha $, violating invariance, so the class $S$ of all singletons that are not elements of their sole element would not qualify for being a set. Thus avoiding the paradox!

\section {A modification of P\'etry Henson Forster schema of invariance}

Here I'll expsoite a modification of Forster's 1995 approach that is already consistent relative to NF, and that gives an equivalent semantic interpretation of stratification. However, the approach is a little bit complex, and it need definition of some notions before explicitly stating it. \bigskip

First we define a canonical indexing $\pi$ on a formula $\phi$ as a function from occurrences of variables in $\phi$ to the naturals, that results from the following indexing procedure: \smallskip

1. If $`x$' is the first occurrence of a variable in $\phi$ [i.e. the most on the left], then   $\pi(`x$')$=|\phi|$, where $|\phi|$ is the number of occurrences of variables in $\phi$

2. If $`x \in y$' is some atomic subformula of $\phi$, then $\pi(`y$'$) = \pi(`x$'$) +1$

3. If $`x = y$' is some atomic subformula of $\phi$, then $\pi(`y$'$) = \pi(`x$'$)$

4. If $\varphi$ is the string of all atomic subformulas of $\phi$ thereof indexed, define $S(\varphi)$ as the first atomic subformula of $\phi - \varphi$ that shares a variable with $\varphi$, or if $\phi-\varphi$ doesn't share any variable with $\varphi$, then $S(\varphi)$ is the first subformula of $\phi-\varphi$ [where: $\phi-\varphi$ is the string that is left from formula $\phi$ after taking $\varphi$ out from it]. Now first we try index the first variable appearing in $S(\varphi)$ by the highest index it received in $\varphi$; if it is not shared with $\varphi$ then we index the second variable of $S(\varphi)$ with the highest index it received in $\varphi$; if both variables of $S(\varphi)$ are not shared with $\varphi$, then we fix the index of the first variable of $S(\varphi)$ to be $|\phi|$.

5- Iterate steps 2,3,4 till the process is exhausted. \bigskip

Define: $i^\pi.rng(x) = |\{\pi(`x$'$): `x$'$ \ occurs \ in \ \phi\}|$, 

in English: the indexing range of a variable $x$ of $\phi$ after $\pi$, is the number of distinct indices of all occurrences of $x$ in $\phi$ according to $\pi$. Now the total sum of all $\pi$-indexing ranges of variables in a formula $\phi$ shall be denoted by $\Sigma \ rng^\pi` \phi$ \bigskip

This way we get $\pi$ to be an indexing function from $occurrences$ of variables in $\phi$ to the $naturals$, that respects requirements of stratification over $occurrences$ of atomic subformulas of $\phi$, and that has the minimal possible $\Sigma \ rng^\pi` \phi$. We call any indexing function that meets these three requirements as $``a \ canonical \ indexing"$ \smallskip

We term $\pi$ as  \emph{\textbf{the}} $canonical \ indexing $ of  $\phi$. \smallskip

Now Quine [9] first defined a stratified formula as a formula whose $variables$ can be indexed by naturals in such a manner as to preserves rules 2,3 depicted above. So stratification is a function on variables, and not just on their occurrences. 

Now $\pi$ would constitute a stratification over formula $\phi$ if and only if $\Sigma \ rng^\pi` \phi = |\{v: v \ is \ variable \ in \ \phi\}|$; that is, the sum indexing of ranges of variables in $\phi$ is equal to the total number of $variables$ in $\phi$. In other words all occurrences of a variable in $\phi$ would be assigned the same index, so $\pi$ would be a function on $variables$ of $\phi$. \smallskip

If $\Sigma \ rng^\pi` \phi > |\{v: v \ is \ variable \ in \ \phi\}|$, then $\pi$ is not a stratified function on $\phi$, and also $\phi$ is not stratifiable! \bigskip

Example: let $\phi$ be $`x \in y \land y \in z \land z \in x$'

The canonical indexing $\pi$ of $\phi$ is: $`x^6 \in y^7 \land y^7 \in z^8 \land z^8 \in x^9$'

So $\Sigma \ rng^\pi` \phi = 2 + 1 +1 = 4 > 3$, so $\phi$ is NOT a stratified formula! \bigskip

Now we define $\phi^f$ as the formula obtained by replacing each occurrence of $`x$'$^n$ of the canonically indexed formula $\phi^\pi$ with $j^n`f(`x$'$)$, where $n$ is the index given to that occurrence of $x$ by the canonical indexing $\pi$ on $\phi$, and $f$ is $2|\phi|$-$setlike$ permutation of $V$. 

Where $f$ is $k$-$setlike$ permutation of $V$ means $f$ is a permutation of $V$ having every $j^n`f$ being a permutation of $V$ also, for all $n \leq k$.

For the sake of short notation we'll say $f$ is $n$-$setlike$ to mean its a permutation over $V$ that is $n$-$setlike$.
Now we come to state the modified invariance schema. \bigskip

\textbf{Modified P\'etry Henson Forster:} if $\phi$ is a formula in which only symbols $`y,u,w$' occur free; and only free, and in which every quantifier is bounded $`\in V$', then:

$\forall u,w \in V:
\forall f (f \ is \ 2|\phi|$-$setlike \Rightarrow \forall y (\phi \longleftrightarrow \phi^f)]) \Rightarrow \exists x \in V( x=\{y: \phi\})$ \bigskip

Now this is a theorem schema of ML which is interepretable in NF. Accordingly NF does prove a semantically motivated re-formulation of it. So a stratified formula is a formula that is invariant under every $2n$-$setlike$ permutation of the universe, where $n$ is the number of occurrences of variables in the formula ($2n \geq$ number of distinct indices used in a stratification).

\section {An aside: Canonical indexing of acyclic formulas}

Earlier works of mine, Holmes and Bowler had shown that stratified comprehension is reducible to acyclic comprehension [2],[3], [4]. An acyclic formula $\phi$ is one in which one can define an acyclic undirected graph $G_\phi$ whose nodes are the $variables$ of $\phi$, and such that for every $occurrence$ of an atomic subformula $`x -- y $' [where $`--$' can be $\in$ or $=$] in $\phi$, there is an edge in $G_\phi$ stretching between the nodes standing for variables $`x,y$'. So it becomes interesting if we can use a version of canonical indexing that can suit defining acyclic formulas in terms of it.

This resembles the general approach above, except that here we don't give separate treatments for identity and membership formulas. We only have $one$ indexing rule of atomic subformulas that is: \smallskip

Rule 1.  for any atomic subformula  $`x \ R \ y$' [where R is either $`\in$' or $`=$'] , if $`x$' occurring in it is labeled first with $n$ then the occurrence of $`y$' in it is indexed with $n+1$, while if the occurrence of  $`y$' in it is indexed first with $k$ then the occurrence of $`x$' in it will be indexed $k+1$. \bigskip

\noindent
The indexing procedure: \smallskip

A. The first variable in the formula receives the natural index 1. \smallskip

B. Now for any initial segment $\varphi$ thereof fully indexed, for $S(\varphi)$ we start indexing the variable in it that has the higher highest index in $\varphi$, if both variables of $S(\varphi)$ have the same highest index in $\varphi$, then we label the first one first, if both variables in $S(\varphi)$ are not shared with $\varphi$, then we'll label the first variable in $S(\varphi)$ with 1. We continue applying rules 1,B until this process is exhausted.\smallskip

Note if a varaible in $S(\varphi)$ is not shared with $\varphi$, then the highest index that it has in $\varphi$ is fixed to 0. \smallskip

Now we follow the same criterion that characterize stratification, so we'd say: \smallskip

$\pi$ is acyclic indexing over formula $\phi$ if and only if 

$\Sigma \ rng^\pi` \phi = |\{v: v \ is \ variable \ in \ \phi\}|$; that is, the sum indexing of ranges of variables in $\phi$ is equal to the total number of $variables$ in $\phi$. In other words all occurrences of a variable in $\phi$ would be assigned the same index, so $\pi$ would be a function over $variables$ of $\phi$. \smallskip

If $\Sigma \ rng^\pi` \phi > |\{v: v \ is \ variable \ in \ \phi\}|$, then $\pi$ is not an acyclic indexing function on $\phi$, even though $\phi$ can be stratified! \smallskip

Example: let $\phi$ be the formula $`x \in y \land z \in y \land k \in x \land k \in z$'

The fully indexed formula $\phi^\pi$ would be: $`x^1 \in y^2 \land z^3 \in y^2 \land k^2 \in x^1 \land k^4 \in z^3$'.

So $\Sigma \ rng^\pi` \phi = 1 + 1 +1+2 = 5 > 4$, so $\phi$ is NOT acyclic formula! Even though it is stratified!\smallskip

This is an aside result of canonical indexing of formulas, I don't know of an approach that can semantically motivate it, like through invariance under permutations or otherwise. It would be interesting to see if there is one along such lines.

\section {References}

\noindent
[1] Al-Johar, Z.A., Short Axiomatization of Stratified Comprehension, pre-print 2020, arXiv:2009.03185v2 [math.LO]

\noindent
[2] Al-Johar, Z.; Holmes M.R., Acyclic Comprehension is equal to Stratified Comprehension, Preprint 2011. http://zaljohar.tripod.com/acycliccomp.pdf

\noindent
[3] Al-Johar, Zuhair; Holmes, M. Randall; and Bowler, Nathan. (2014). "The Axiom Scheme of Acyclic Comprehension". Notre Dame Journal of Formal Logic, 55(1), 11-24. http://dx.doi.org/10.1215/00294527-2377851

\noindent
[4] Al-Johar, Zuhair. Stratified formulas are equivalent to Acyclic formulas: A review. pre-print 2020, arXiv:2009.13274v2 [math.LO]

\noindent
[5]Ehrenfeucht, A., and A. Mostowksi, “Models of axiomatic theories admitting automorphisms,” Fundamenta Mathematicae, vol. 43 (1956), pp. 50–68. Zbl 0073.00704. MR 0084456. 574

\noindent
[6] Forster, T. E., Set theory with a universal set: exploring an untyped uni-verse (2nd ed.), Oxford Logic Guides no. 31, Clarendon Press, Oxford,1995.

\noindent
[7] Holmes, M. R., ``Symmetry as a criterion for comprehension motivating Quines `New Foundations'", Studia Logica, vol. 88, no. 2 (March 2008).

\noindent
[8] Holmes, M.R., Symmetric comprehension revisited, pre-print 2020.

\noindent
[9] Quine, W. v. O., ``New foundations for mathematical logic", American Mathematical Monthly, vol. 44 (1937), pp. 70-80.

\noindent
[10]Quine, Willard Van Orman (1951), Mathematical logic (Revised ed.), Cambridge, Mass.: Harvard University Press, ISBN 0-674-55451-5, MR 0045661

\noindent
[11]Rosser, Barkley (1942), "The Burali-Forti paradox", Journal of Symbolic Logic, 7 (1): 1–17, doi:10.2307/2267550, JSTOR 2267550, MR 0006327

\noindent
[12]Wang, Hao (1950), "A formal system of logic", Journal of Symbolic Logic, 15 (1): 25–32, doi:10.2307/2268438, JSTOR 2268438, MR 0034733

\noindent
[13]Whitehead, Alfred North and Bertrand Russell (1963). Principia Mathematica. Cambridge: Cambridge University Press.

\end{document}